\documentclass[11pt]{article}
\usepackage{amsmath}
\usepackage{amsfonts}
\usepackage{amsthm}
\usepackage{verbatim}

\begin{document}

\def\N{\mathbb{N}}
\def\F{\mathbb{F}}
\def\Z{\mathbb{Z}}
\def\R{\mathbb{R}}
\def\Q{\mathbb{Q}}
\def\H{\mathcal{H}}

\parindent= 3.em \parskip=5pt

\centerline{\bf{A short history of some recent research }} 
\centerline{\bf{ on continued fractions in function fields}}
\vskip 0.5 cm
\centerline{\bf{by A. Lasjaunias }}

\vskip 0.5 cm {\bf{Abstract.}} The goal of this survey paper is to present, in chronological order, certain research works on continued fractions in power series fields over a finite field, all of them being derivated from some examples introduced thirty years ago by Mills and Robbins.
\vskip 0.5 cm
 For general information concerning the frame of this note, diophantine approximation and continued fractions in power series fields, and also more references, the reader may consult \cite{L}, \cite{S}  and \cite[Chap.~9]{T}.

In the sequel, $\F_q$ is the finite field of characteristic $p$, with $q$ elements. We denote $\F(q)$ the field of power series in $1/T$ over $\F_q$, where $T$ is a formal indeterminate. In this field every element $\alpha$ is expanded as a finite (if $\alpha$ is rational, i.e. $\alpha \in \F_q(T)$) or infinite continued fraction $\alpha=[a_0,a_1,\cdots,a_n,\cdots]$, where the partial quotients $a_n\in \F_q[T]$ are nonconstant polynomials except possibly for $a_0$. Moreover if $\alpha$ is an irrational element in $\F(q)$ satisfying an algebraic equation of the form $A\alpha^{r+1}+B\alpha^{r}+C\alpha+D=0$, where $A,B,C$ and $D$ are in $\F_q[T]$ and $r$ is a power of the characteristic $p$, then we say that $\alpha$ is hyperquadratic (of order $t$ if $r=p^t$). 

The starting point of the study of continued fractions for algebraic power series over a finite field is certainly due to Baum and Sweet \cite{BS1,BS2}, in the midle of the 1970's. Working in the same research institute (IDA, Princeton), Mills and Robbins follow this line by publishing another very important paper \cite{MR}. This last paper contains a study of a famous example introduced by Baum and Sweet, but also several other interesting and somehow mysterious continued fraction expansions. In \cite[p.~399-402]{MR}, they introduce particular algebraic continued fractions in $\F(p)$ having all partial quotients of degree one. For $p\geq 5$, these examples are such that $a_n=\lambda_nT$ for $n\geq 1$, where $\lambda_n\in\F_p^*$. Also for $p=3$, a particular example is given where $a_n=\lambda_nT+\mu_n$ for $n\geq 1$, with $\mu_n\in\F_3$. Note that these algebraic elements are all hyperquadratic and this was used to obtain their explicit continued fraction expansion. Finally in the same article \cite[p.~403-404]{MR}, they make the following observation: the equation $x^4+x^2-Tx+1=0$ has a solution in $\F(p)$ for each $p$ and the continued fraction expansion for this solution appears to be remarkable for two values of $p$: if $p=3$ and $p=13$.

In another article \cite{BR}, about ten years after, Buck and Robbins turned to this strange quartic equation introduced at the end of \cite{MR}. First, they prove the conjecture, made earlier, on the explicit continued fraction in the case $p=3$. Secondly, they make a full conjecture (based on computer observation) on the case $p=13$.

In \cite{L1}, just shortly after the publication of Buck and Robbins, I did study the same continued fractions for the solution of this quartic equation. In the case $p=3$, I could give an alternative proof of the conjecture, giving rise to a posible generalization. Moreover I made the observation that the root of the quartic equation is not hyperquadratic if $p=3$, while it is so if $p=13$. Actually in my PhD thesis (1996), I could prove that the root, in $\F(p)$ for $p>3$, of $x^4+x^2-Tx-1/12=0$ is  hyperquadratic of order 1 ($r=p$) if $p\equiv 1 \mod 3$ and hyperquadratic of order 2 ($r=p^2$) if $p\equiv 2 \mod 3$. Due to the use of computer calculations, this was only true with a limitation on the prime $p$. 

In \cite{L2}, I was inspired by the example in $\F(3)$, introduced in \cite{MR}, having all partial quotients of degree one. I could describe a large family of hyperquadratic elements also in $\F(3)$ having the same property. The method used in \cite{L2}, was later extended in joint works with J-J. Ruch \cite{LR1,LR2}, trying to understand the existence of these algebraic continued fractions with partial quotients of degree one in the larger context of $\F(q)$. Unfortunately, we could only obtain partial results. So I decided to abandon this way.

A few years later, I turned to the question of deciding if a quartic power series is hyperquadratic or not (Remark that quadratic and cubic elements are always hyperquradratic). Then I contacted A. Bluher, from the NSA, and in a joint work we could obtain a general result \cite{BL}, removing the limitation on the size of $p$ in the result introduced in my thesis.

I was trying to understand the origin of the continued fraction for the quartic in the case $p=13$. Generalizing the problem to the equation $x^4+x^2-Tx-1/12=0$ for various  $p\equiv 1 \mod 3$, I noticed that these continued fraction expansion were all built from a particular polynomial of the form  $P_k=(T^2+u)^k$, where $u\in \F_p^*$ and $k\in\N$. Hence, the solution of the Mills-Robbins quartic for $p=13$ was based on $P=(T^2+8)^4$. Moreover, I noticed that the elements mentioned above and introduced by Mills and Robbins, having partial quotients of degree one, were also built from such a polynomial; in this case $P=(T^2+4)^{(p-1)/2}$. These expansions were introduced in \cite{L3}, and I will call them now $P_k-expansions$.

However, the study begun in \cite{L3} was not sufficient to obtain a proof of Buck and Robbins conjecture for the solution of the quartic with $p=13$. This was only obtained in \cite{L4} by deepening the method. Some particular continued fractions, built from the polynomials $P_k$, are such that the partial quotients are all proportional to some polynomials derived from a particularly simple sequence $(A_n)_{n\geq 0}$. Therefore, they are now called $P_k-expansions$ $of$ $type$ $A$. Finally a more general description of these $P_k-expansions$, in the wider context of $\F(q)$, was presented in \cite{L5}. It happens that, for a particular choice of the parameter $k$, this sequence $(A_n)_{n\geq 0}$ is constant, with $A_n=T$, and this leads to the original examples described in \cite{MR}, with partial quotients of degree one.

In 2009, I met D. Gomez at a conference in Dublin. We had a correspondance and the goal was to describe the continued fraction for the solution of the quartic equation  
 $x^4+x^2-Tx-1/12=0$ for all $p>3$. The case $p=13$ was known. The cases $p\equiv 2 \mod 3$ were still mysterious. Finally, I published a short note \cite{L6}, resuming the cases $p\equiv 1 \mod 3$, with a proved description of this expansion, again with a limitation on the size of $p$. Meanwhile, we worked with Gomez on a generalisation of hyperquadratic continued fractions with partial quotients of degree one in $\F(3)$. These are also particular $P_k-expansions$ $of$ $type$ $A$. A note \cite{GL} was published.

The existence of these algebraic continued fractions with partial quotients of degree one was a remarkable result in \cite{MR}. Extending the work with Gomez \cite{GL} to an arbitrary finite base field, I realised that the general case, in odd characteristic, could probably be treated in this way, getting over the unfructuous attemps made with Ruch several years earlier. I had discussions on this matter wir A. Firicel during some months. The problem was quite technical. In february 2014, J-Y. Yao was visiting Bordeaux university and we took this oportunity to work on this problem. Finally, we could conclude and a joint note was published \cite{LY}.

During the summer 2014, I turned to the generalization of Robbins'quartic. In october 2013, Kh. Ayadi was in Bordeaux for a conference and I proposed him to work on this matter. We used extensive computer calculations and we could finally describe the pattern of the continued fraction for the solution in the second case $p\equiv 2 \mod 3$. In this case, the expansion is also built from a polynomial $P_k$ as above. However, the partial quotients are all proportional to some polynomials derived from a second sequence $(B_n)_{n\geq 0}$ and the element is hyperquadratic of order 2 ($r=p^2$). Therefore, we could define $P_k-expansions$ $of$ $type$ $B$ of second order. A note, explaining this but containing mostly conjectures, was submitted at the end of 2014. It was published recently \cite{AL}. 

Shortly after Mills and Robbins'publication, J.-P. Allouche \cite{JPA} could show, for the continued fractions with partial quotients of degree one introduced there, that the sequence of partial quotients is automatic. During the summer 2015, I thought about the possible connection between certain algebraic continued fractions and automatic sequences in a finite field. I thought the automaticity of the sequence of the leading coefficients of the partial quotients could depend on the caracter hyperquadratic or not of the continued fraction. We had an example non-hyperquadratic with the solution of Robbins'quartic for $p=3$ and several examples hyperquadratic (besides Allouche's result cited above, other simple examples in characteristic 2 were known but unpublished). J.-Y. Yao being a specialist on the matter of automatic sequences, we decided to work together on this subject. We could first prove that the sequence derived from the solution of Robbins'quartic for $p=3$, which has only two values 1 or 2, is non-automatic. Then, we submitted a joint work \cite{LY1} where several examples of automatic sequences linked to hyperquadratic continued fractions were presented, as well as this counter-example in $\F_3$. Note that this two-valued sequence was later on studied more deeply, leading to a nice continued fraction in characteristic zero \cite{ABL}. At the end of the year 2015, I also realised that the technique introduced to treat hyperquadratic continued fractions with partial quotients of degree one in odd characteristic could also be applied to the characteristic 2. This underlined the role played in both cases by the analogue of the golden mean in the formal case. A note was submitted \cite{L7} and this note was the starting point of another joint work with Yao \cite{LY2}, where the link between automaticity and certain continued fractions over a finite field was studied. In this last note, general results on recurrent and automatic sequences were presented. But we also could obtain the automaticity of all the sequences derived from $P_k-expansions$ $of$ $type$ $A$, in the more general setting, including the first cases introduced in \cite{MR} and studied by Allouche in \cite{JPA}.

\vskip 0.5 cm
\begin{tabular}{ll}Alain LASJAUNIAS\\Institut de Math\'ematiques de Bordeaux  CNRS-UMR 5251
\\Universit\'e de Bordeaux \\Talence 33405, France \\E-mail: Alain.Lasjaunias@math.u-bordeaux.fr\\\end{tabular}

\end{document}